# THE TARGET-MATRIX OPTIMIZATION PARADIGM FOR HIGH-ORDER MESHES*

VESELIN DOBREV†, PATRICK KNUPP‡, TZANIO KOLEV†, KETAN MITTAL§, AND VLADIMIR TOMOV†

**Abstract.** We describe a framework for controlling and improving the quality of high-order finite element meshes based on extensions of the Target-Matrix Optimization Paradigm (TMOP) of [24]. This approach allows high-order applications to have a very precise control over local mesh quality, while still improving the mesh globally. We address the adaption of various TMOP components to the settings of general isoparametric element mappings, including the mesh quality metric in 2D and 3D, the selection of sample points and the solution of the resulting mesh optimization problem. We also investigate additional practical concerns, such as tangential relaxation and restricting the deviation from the original mesh. The benefits of the new high-order TMOP algorithms are illustrated on a number of test problems and examples from a high-order arbitrary Eulerian-Lagrangian (ALE) application [6]. Our implementation is freely available in an open-source library form [31].

**Key words.** target-matrix optimization paradigm, high-order meshes, high-order finite elements, mesh optimization.

**1. Introduction.** High-order methods are becoming increasingly important in computational science due to their potential for better simulation accuracy and favorable scaling on modern architectures [37, 13, 12, 9, 10]. A vital component of such methods is the use of high-order representation not just for the *physics* fields, but also for the geometry, represented by a high-order computational mesh. High-order finite element meshes in particular, can be very beneficial in a wide range of applications, where e.g. radial symmetry preservation, or alignment with physics flow or curved model boundary is important [36, 14, 30, 20, 7]. Such applications can utilize *static* meshes, where a good-quality high-order mesh needs to be generated only as an input to the simulation, or *dynamic* meshes, where the mesh evolves with the problem and its quality needs to be constantly controlled. In both cases, the quality of high-order meshes can be difficult to control, because their properties vary in space on sub-zonal level. Such control is critical in practice, as poor mesh quality leads to small time steps or simulation failures.

In this paper we propose theory that defines high-order mesh quality, and using that theory we produce general and flexible mesh optimization tools for current and future high-order applications, as well as low-order applications that use high-order meshes. Our work is based on an extension of the Target-Matrix Optimization Paradigm (TMOP) of [24], incorporating also some ideas from variational minimization (VM) [8, 38, 16]. This approach provides a framework that can be tailored to the application's needs and is easy to use and freely available in a library form [31]. The methods we develop are general (applicable to any element obtained through mapping from a *reference* element), but we illustrate them on quad/hex meshes mostly focusing on purely geometric mesh optimization.

There is extensive literature on the generation, optimization, and adaptation of unstructured meshes containing first-order element types. In comparison, there is a paucity of literature on generating meshes containing high-order elements; some examples are [39, 35, 21, 26, 28, 33, 19, 25, 34]. While our focus on is on purely geometric optimization, we note that there is also relatively little published concerning solution-adapted meshes containing high-order elements [40, 30, 29, 1], not all of which employ $r$-adaptivity.

† Center for Applied Scientific Computing (CASC), Lawrence Livermore National Laboratory, Livermore, CA 94551, dobrev1@llnl.gov, tzanio@llnl.gov, tomov2@llnl.gov
‡ Dihedral LLC, zg37rd@gmail.com
§ University of Illinois Urbana-Champaign, kmittal2@illinois.edu

TMOP is a general approach for controlling mesh quality, where mesh nodes (vertices in the low-order case) are moved so-as to optimize a multi-variable objective function that quantifies global mesh quality. Quality, in turn, is measured at a collection of points within elements of the mesh by means of a local quality metric. The metric is a function of a pair of matrices representing the Jacobian of the element map and the ideal Jacobian (or target). The Jacobian matrix is a function of the mesh node coordinates. By selecting appropriate target-matrices and local quality metrics, optimization improves mesh quality in specific ways. TMOP is distinguished from similar methods by its emphasis on target-matrix construction methods that permit a greater degree of control over the optimized mesh. Element quality is not directly computed but is itself an average over quality at local points within the element. It is this property which allows TMOP to be extended to the optimization of high-order meshes.

Controlling the quality of high-order meshes requires few adaptations and changes to the TMOP theory developed for linear elements. However, TMOP for high-order meshes has not been seriously implemented in code previously and thus its actual performance in this context was unknown prior to this work. In this paper we explore the application of TMOP to high-order meshes and report on the performance of the resulting high-order TMOP algorithms for meshes arising in a high-order arbitrary Eulerian-Lagrangian (ALE) application [6].

Our approach is based on a combination of TMOP and VM. In variational minimization, the mesh optimization problem is written as a sum of integrals over the elements with the quality metric as the integrand. Consequently, the local optimization problems are coupled through the high-order nodes of the mesh shared between neighboring elements. An optimal target element is specified for each element, and the local quality is measured based on the mapping to that element. Two natural choices for the target elements are (a) a single optimal element, e.g., the unit square, an equilateral triangle, the unit cube, etc. and (b) the elements of a *good* reference mesh that has the same connectivity as the mesh to be optimized, with the goal being to map the *good* mesh to a domain with a different shape or to $r$-adapt the mesh. The integrals in the optimization problem are discretized using a quadrature rule, and the quadrature points play the role of (weighted) sample points.

As with many methods for mesh generation and mesh quality improvement, our approach may require the user to select values of certain numerical parameters which appear in the target-matrices, the composite quality metrics, or in the number of sample points. Some degree of trial-and-error on the part of the user may be necessary to determine good values for the parameters. To make this selection process easier, the parameters within the target-matrices and composite metrics have a clear meaning and are often scaled to be in the range $[0, 1]$. Although our method is thus not always fully automatic, these user-parameters lend considerable flexibility and capability to the approach.

Section 2 reviews the representation of the high-order mesh, basic TMOP components, and the objective function. Section 3 describes the TMOP target-matrices, quality metrics, quadrature point selection, and how we numerically solve the opti-mization problem. Section 4 considers tangential relaxation of boundary node positions, composite metrics, and methods for restricted node movement. Section 5 presents numerical results which demonstrate the ability to control high-order meshes using basic TMOP target-matrices and quality metrics.

**2. Preliminaries.** The goal of this section is to introduce the discrete mesh representation that is used throughout the paper, the main components of TMOP, and the corresponding notation. All these concepts are already described in previous papers, but we repeat them here for completeness. The changes and adaptations in the TMOP components for high-order meshes are discussed in the following section.

**2.1. Discrete representation of the high-order mesh.** Let $d \in \{1, 2, 3\}$ be the space dimension and $\mathcal{V} \subset [H^1(\Omega)]^d$ be a space of continuous finite element functions defined on $\Omega$. In particular,

$$\mathcal{V} = \begin{cases} \{v \in [C^0(\Omega)]^d \mid v_{|E} \in (Q_k)^d\} & \text{for quadrilateral/hexahedral meshes}, \\ \{v \in [C^0(\Omega)]^d \mid v_{|E} \in (P_k)^d\} & \text{for triangular/tetrahedral meshes}, \end{cases}$$

where $Q_k$ is the space of polynomials of degree at most $k$ in each variable, and $P_k$ is the space of polynomials of total degree at most $k$. We discretize mesh positions by the high-order finite element function $x \in \mathcal{V}$, i.e., on each element $E$ we use the expansion

$$x(\bar{x}) = \sum_{i=1}^{N} \mathbf{x}_i w_i(\bar{x}), \tag{2.1}$$

where $\{w_i\}_{i=1}^{N}$ is the basis of $\mathcal{V}$ on $E$, $\bar{x}$ are the fixed coordinates of the reference element, and $\mathbf{x}$ is the coordinate vector, which contains the values of the finite element degrees of freedom, of size $Nd$ (each $\mathbf{x}_i$ in (2.1) is of size $d$). The $d \times d$ Jacobian of the transformation from reference to physical coordinates reads

$$A(x) = \frac{\partial x}{\partial \bar{x}} = \sum_{i=1}^{N} \mathbf{x}_i \nabla w_i(\bar{x}). \tag{2.2}$$

**2.2. Basic TMOP components.** Assuming that all reference → physical mappings are differentiable so that (2.2) is well defined, the first component of TMOP is a set of *sample points* within the reference element. These points are extended to the physical space by the (2.1) mapping. At each sample point (inside each mesh element), TMOP uses two Jacobian matrices:

- The Jacobian matrix $A_{d \times d}$ of the transformation from reference to physical coordinates, given by (2.2).
- The *target matrix*, $W_{d \times d}$, which is the Jacobian of the transformation from reference to the coordinates of a *target element*. The target elements are defined according to a user-specified method prior to the optimization; they define the desired elements in the optimal mesh.

The set of target matrices remains fixed during the optimization procedure while $A$ can change since it depends on the nodal coordinates. The *weighted Jacobian* matrix, $T_{d \times d}$, defined by $T(x) = A(x)W^{-1}$, represents the Jacobian of the transformation between the target and the physical coordinates. This matrix is used to define the *local quality measure*, $\mu(T)$, see Section 3.2. The quality measure can evaluate shape, size, or alignment of the region around the sample point of interest.

The goal of TMOP is to minimize a global *objective function*, $F(\mu)$ (see Section 2.3), which combines the quality measures in the chosen set of sample points throughout the mesh, so that the final mesh elements are as close as possible to the shape, size, or alignment of their targets. The original mesh topology is preserved, i.e., optimization is performed via node movement. The optimization procedure is also constrained by boundary conditions. Additional details about TMOP and properties of the local quality metrics can be found in [22, 23, 24].

**2.3. Objective functions.** The objective function (or norm) $F(\mu)$, in which the mesh is optimized, can be defined in various ways. Two examples of such objective functions are:

- The *pointwise* objective function, used by the original TMOP formulation:

$$F_P(x) := \sum_{E \in \mathcal{E}} \sum_{x_s \in S_E} c_s^E \, \mu(T(x_s)) \,, \tag{2.3}$$

  where $\mathcal{E}$ is the set of high-order elements in the mesh, $S_E$ is the set of sample points within element $E$, $c_s^E$ is a sample point *trade-off coefficient*, and $T(x_s)$ is evaluated at the sample point $x_s$ of element $E$.

- The *variational* objective function:

$$F_V(x) := \sum_{E \in \mathcal{E}} \int_{E_t} \mu(T(x)) dx_t = \sum_{E \in \mathcal{E}} \sum_{x_q \in Q_E} w_q \, \det(W(x_q)) \, \mu(T(x_q)) \,, \tag{2.4}$$

  where $E_t$ is the target element corresponding to the physical element $E$, $Q_E$ is the set of quadrature points for element $E$, $w_q$ are the corresponding quadrature weights, and both $T(x_q)$ and $W(x_q)$ are evaluated at the quadrature point $x_q$ of element $E$.

Optimal mesh configuration is achieved by minimizing the objective function. The two formulations are equivalent when the sample and quadrature points are chosen to have the same positions, and $c_s = w_q \det(W(x_q))$ for the corresponding $s$ and $q$ indices. Both formulations admit the inclusion of space-dependent weights for each quadrature point, which would make some parts of the domain more important than others. In this paper, we focus on the $F_V(x)$ objective function.

The existence of a minimum for (2.4) has been explored theoretically in the context of hyperelasticity [5], and applied to variational mesh optimization in [17, 18]. The cited work shows that there exists an elastic deformation (i.e., mesh configuration $x$) that minimizes the deformation energy functional $(F)$, provided the elastic potential $(\mu)$ is a polyconvex function and satisfies certain growth conditions.

**3. High-order mesh optimization with TMOP.** This section discusses the application of TMOP to high-order meshes and highlights the algorithmic components that are specifically relevant in this context.

**3.1. TMOP target matrices.** We start by extending the notion of target matrices to sub-element level. In particular, we are free to specify the target matrix $W$ at every quadrature point, provided $\omega = \det(W) > 0$. Defining target matrices at quadrature points provides a way to enforce sub-element features in the high-order mesh elements. As an (artificial) example, consider a high-order mesh that contains a single element. Suppose the goal of the optimization is to stack the element's quadrature points towards the middle of the element, while the boundary of the element stays fixed. Clearly, the target matrices of the inner-most quadrature points must incorporate the notion of smaller size. We present a similar (and more practical) example in Section 5.2.4.

In the left-hand side of (2.4), we assumed that we are given a target element, $E_t$, for every element $E$ in the mesh. Then, based on the mapping from the reference element that defines $E_t$, we computed the target Jacobian matrices $\{W(x_q)\}$ at the quadrature points. We can generalize this definition of the target Jacobians, $\{W(x_q)\}$, by allowing them to be defined directly. In other words, there is no need to assume that $\{W(x_q)\}$ are obtained by differentiating a given target element mapping. Thus, we can take the right-hand side of (2.4) as a more general definition of the variational objective function, where $\{W(x_q)\}$ are given as input, e.g., constructed using a *target Jacobian construction* algorithm. In this setting, the weighted Jacobian matrices, $T$, are computed as $T = AW^{-1}$.

Several target sets are used in this paper, depending on the application. Below we list their 2D versions, but all targets are easily extendable to 3D.

The first target-matrix set corresponds to the shape of an isotropic quadrilateral element, namely a unit square. The resulting set of target matrices are $W_1 = I$, the identity matrix, constant over all quadrature points.

The second target-matrix set corresponds to the shape of an isotropic quadrilateral element having a prescribed area. In particular, the target in this case is the constant $W_2 = \sqrt{\bar{\alpha}}\, I$, where

$$\bar{\alpha} = \frac{1}{N_{\mathcal{E}}} \int_\Omega 1 \ dx \,, \tag{3.1}$$

with $N_{\mathcal{E}}$ being the number of elements in the mesh. The intent of this constant target is to produce equal-area, square elements. In this case, if the optimization is perfect, i.e., $T = I$ everywhere, we would have $\det(A) = \det(TW) = \bar{\alpha}$. Note that (3.1) can be scaled additionally by a space-dependent factor. Such targets are useful when trying to enforce certain sizes in specific regions of the mesh, see Section 5.2.4.

The third 2D target-matrix set corresponds to the shape of an isotropic quadrilateral element having a local area the same as the corresponding local area in the initial mesh (assuming the local areas are non-negative). The target in this case is the non-constant $W_3(x(x_0)) = \sqrt{\det(A(x_0))}\, I$, where the determinant is evaluated on the initial mesh. The intent of this target is to produce square elements while preserving the local sizes of the initial mesh.

## 3.2. TMOP quality metrics in the high-order case.
A main component in TMOP is the local quality metric. A TMOP local quality metric is a functional $\mu = \mu(T)$ from $\mathcal{D} \subseteq M_d$ to the set of non-negative numbers. Here $M_d$ is the set of all real $d \times d$ matrices $d = 2, 3$ [1]. For high-order meshes, the matrix $T$ generally varies within each element, hence it is important to evaluate $\mu(\cdot)$ at various positions (quadrature points) within the element. We say that a quality metric is well-behaved when $\mu(T) \geq 0$, with $\mu(T) = 0$ if and only if $T = T_m$, where $T_m$ belongs to some prescribed minimizing set $\mathcal{M} \subset \mathcal{D}$.

Depending on $\mathcal{M}$, TMOP metrics can be classified into a variety of canonical *types*. The different metric types are defined through the standard decomposition of the target Jacobian matrix into four $d \times d$ matrices [22]:

$$W = [\text{volume}][\text{orientation}][\text{skew}][\text{aspect ratio}] \,. \tag{3.2}$$

---

[1] Extensions to manifolds can also be made, but we are not concerned with that topic in this paper.

A metric $\mu(T)$ is defined in a way that measures the difference between $A$ and $W$ only in particular components of (3.2), and is invariant of the others. The metric type depends on the set of chosen components.

- *Shape (Sh)* metrics control [skew] and [aspect ratio]. They are minimized when $A$ is a scaled rotation of $W$. Shape metrics are invariant of [rotation] and [volume] of $A$. Example of a shape metric is $\mu_2 = \frac{|T|^2}{2\tau} - 1$, where $\tau = \det(T)$, and $|T|^2 = \operatorname{tr}(T^tT)$.

- *Size (Sz)* metrics control [volume]. They are minimized then $\det(A) = \det W$, and are invariant to the other components of the decomposition. Example of a size metric is $\mu_{77} = 0.5\left(\tau - \frac{1}{\tau}\right)^2$.

- *Alignment (Al)* metrics control [orientation] and [skew]. They are minimized when $A = WD$, where $D$ is a diagonal matrix. Example of an alignment metric is $\mu_{30}(A, W) = |\boldsymbol{a_1}||\boldsymbol{w_1}| - (\boldsymbol{a_1} \cdot \boldsymbol{w_1}) + |\boldsymbol{a_2}||\boldsymbol{w_2}| - (\boldsymbol{a_2} \cdot \boldsymbol{w_2})$, where $\boldsymbol{a}$ and $\boldsymbol{w}$ are the column vectors of $A$ and $W$.

- The above three major types can be combined, e.g., into *Shape+Size (SS)* metrics like $\mu_7 = |T - T^{-t}|^2$, or *Shape+Size+Alignment (SSA)* metrics like $\mu_{14}(T) = |T - I|^2$. SS metrics are minimized when $A$ is a rotation of $W$, while SSA metrics are minimized only when $A = W$, i.e., $T = I$.

Presenting a full list of our metrics and comparisons of their properties (e.g., convexity and polyconvexity) will be the subject of a future work.

In addition to the metric type, another useful classification pertains to metric *category*. Metric category has to do with the domain of metric and the method it uses to avoid the creation of inverted optimal meshes. The highest level metric categories are (a) *non-barrier* metrics and (b) *barrier* metrics. Under the non-barrier category is a specialized category known as the *pseudo-barrier* category of metrics. Under the barrier category are the *ideal* and *shifted-barrier* categories. Metrics whose category is ideal-barrier have domain $\mathcal{D}$ equal to the set of all $d \times d$ matrices whose determinant is positive. From a theoretical point of view, the ideal barrier category is central, while from a practical viewpoint, the shifted- and pseudo-barrier categories are important. Since none of the initial meshes used in this study contain negative area at the quadrature points, we are concerned exclusively with ideal barrier metrics in the paper.

Table 3.2 summarizes specific metrics used in this study to produce numerical results in Section 5. In the formulas defining the local quality metrics $\tau$ is the determinant of T, $\det(T)$, and $|T|$ is the Frobenius norm of $T$.

Note that metric $\mu_1$ does not fit the formal definition of any of the above types. This metric optimizes towards the target shape, however, it is generally not well behaved, as $\mu_1(T)$ would always approach zero when the size of the physical element decreases, independent of the used target. Nevertheless, in practice the minimal size of elements is usually restricted, e.g., because of boundary conditions. Thus, one can think of $\mu_1$ as *Shape+Smallest Possible Size* metric. We include $\mu_1$ in this study as its use has proven beneficial in cases that include additional limitations about the decrease of the local mesh sizes, for example in Sections 5.2.3 and 5.2.4.

**3.3. TMOP quadrature point selection in the high-order case.** As previously noted, the local quality metrics are evaluated at quadrature points within a given element. In general, one can use any quadrature rule for integrating (2.4). The quadrature rule is selected prior to optimization. For example, for a 2D target element, one could select the tensor product of $N$ uniformly distributed points in each of two directions, giving a total of $N^2$ quadrature points. Alternatively, one can

| ID | Dimension | Type | Category | Formula |
|----|-----------|------|----------|---------|
| 1 | 2 | none | No Barrier | $\mu_1 = \mid T \mid^2$ |
| 2 | 2 | Sh | Ideal Barrier | $\mu_2 = \frac{|T|^2}{2\,\tau} - 1$ |
| 9 | 2 | SS | Ideal Barrier | $\mu_9 = \tau \mid T - T^{-t} \mid^2$ |
| 77 | 2 | Sz | Ideal Barrier | $\mu_{77} = \frac{1}{2}(\tau - \frac{1}{\tau})^2$ |
| 303 | 3 | Sh | Ideal Barrier | $\mu_{303} = \frac{|T|^2}{3\,\tau^{2/3}} - 1$ |

TABLE 3.1

*Examples of specific quality metrics of different type, dimension and category. We denote: $\tau = \det(T)$, and $|T|^2 = \operatorname{tr}(T^t T)$.*

use any of the standard Gauss-Lobatto or Gauss-Legendre quadrature rules. In our experience, we have not observed any systematic differences between using uniform, Gauss-Lobatto, or Gauss-Legendre spacing and weights with equal number of points.

It is critical with high-order meshes to choose $N$ sufficiently large in order to obtain a good numerical solution to the optimization problem. At a minimum, the number of quadrature points $N$ should be one more than the order of the element map. Experience shows that this generally suffices for map orders 1-3. For higher order maps the number of quadrature points may need to be even larger in order to avoid poor mesh quality resulting from under-sampling of the quality within the element. If too few sample points are used, the numerical algorithm can produce inverted and/or non-smooth optimized meshes, even when the underlying quality metric is sound. On the other hand, choosing $N$ too large gives diminishing returns in terms of quality improvement and increases overall computational expense. Figure 3.1 shows that utilizing an insufficient number of quadrature points (in this case $N = 5$ points in each direction) can lead to inverted elements when smoothing a fourth order mesh. This is avoided when $N = 6$ is used. In both cases, the optimization algorithm converges successfully. Both results use the Gauss-Lobatto spacing and quadrature coefficients.

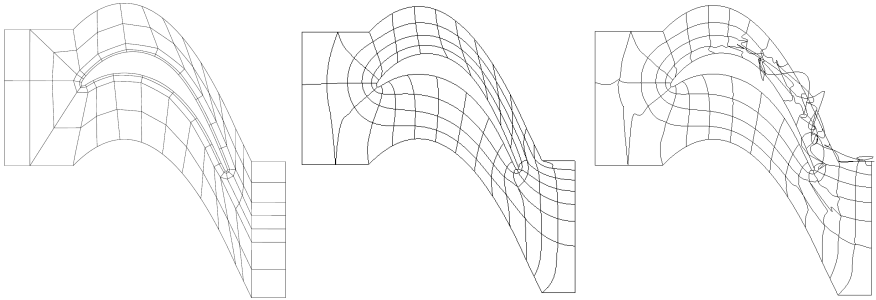

FIG. 3.1. *Comparison of the original mesh (left) with smoothed meshes: $N = 6$ (center) and $N = 5$ (right).*

**3.4. Solving the TMOP optimization problem.** The global TMOP quality functional defined in the previous sections is a nonlinear function with respect to the node positions of the high-order mesh. Thus, to improve the mesh quality, we need to apply methods from the field of nonlinear optimization. Full exploration and evaluation of this field is outside the scope of this article. Instead, we describe two simple optimization strategies that work well in practice and were used in our numerical tests.

For both *pointwise* (2.3) and *variational* (2.4) objective functions, one can apply a derivative-free optimization method. Our derivative-free implementation iterates over all nodes and consecutively displaces each node $\mathbf{x}_i$ to a position that minimizes the local patch mesh quality, $Q(\mathbf{x}_i)$, around the node. For (2.3) and (2.4), respectively, we can use

$$Q(\mathbf{x}_i) = \sum_{E \in N_i} \sum_{\mathbf{x}_s \in E} c_s^E \, \mu(T_s^E), \quad \text{and} \quad Q(\mathbf{x}_i) = \sum_{E \in N_i} \int_{E_t} \mu(T(x)) dx_t, \qquad (3.3)$$

where $N_i$ is the set of elements (neighbors) that contain the node $\mathbf{x}_i$. The local minimization of $\mathbf{x}_i$'s coordinates is performed by the Nelder-Mead simplex method [32, 27]. The outer iterations over the mesh nodes are terminated when

$$\max_k \left| \mathbf{x}_k^{m+1} - \mathbf{x}_k^m \right| < t,$$

where $t > 0$ is a user-specified tolerance that depends on the volume of the domain (we use $t = 10^{-5}$ for domains of unit volume), and $m$ indicates the $m$-th iteration of the procedure. The main advantages of this derivative-free method are its simplicity and the ability to iterate over subsets of the domain. Its disadvantages are slow convergence and difficult parallel implementation, as each displacement in the outer iteration depends on the previously performed displacements (essentially a Gauss-Seidel type of iteration). Optimizing this method will be an object of future work.

As an alternative, we utilize Newton's method to solve the critical point equations, $\partial F(\mathbf{x})/\partial \mathbf{x} = 0$. For our variational objective function (2.4), this derivative has the form

$$\frac{\partial F_V(\mathbf{x})}{\partial \mathbf{x}_i} = \sum_E \int_{E_t} \left( \frac{\partial A(\mathbf{x})}{\partial \mathbf{x}_i} W^{-1} \right) : \frac{\partial \mu(T)}{\partial T} dx_t = \sum_E \int_{E_t} \left( \nabla \tilde{\phi}_j W^{-1} \right) : \frac{\partial \mu(T)}{\partial T} dx_t,$$

where we used the identities

$$A(\mathbf{x}) = \sum_j \mathbf{x}_j \nabla \tilde{\phi}_j(\tilde{x}), \quad \left[ \frac{\partial \mu}{\partial T} \right]_{kl} = \frac{\partial \mu}{\partial T_{kl}}.$$

An expression for the Hessian of $F_V$, which is required for the computation of the system Jacobian, can be derived similarly, involving the computation of the second derivatives of $\mu(T)$ with respect to $T$. For easier derivations of these derivatives, it is beneficial to define the quality metric $\mu(\cdot)$ in terms of the 2D or 3D matrix invariants (or other quantities with known derivatives). All results presented in Section 5 minimize the *variational* objective (2.4) with this Newton-based method.

Some metrics require $\det(T) > 0$ at every sample or integration point, which constrains the admissible configurations of the mesh nodes. In our derivative-free method, this is addressed by setting $\mu(T)$ equal to a big value whenever $\det(T) \le 0$, causing the local Nelder-Mead method to avoid such configurations. When Newton's method is used, we consider the following scaling of the update:

$$\mathbf{x}^{n+1} = \mathbf{x}^n - \alpha \left[ \frac{\partial^2 F_V(\mathbf{x}^n)}{\partial \mathbf{x}^2} \right]^{-1} \frac{\partial F_V(\mathbf{x}^n)}{\partial \mathbf{x}},$$

and decrease the scaling factor $\alpha$ until the new configuration The $\alpha$ factor is also decreased when $\mathbf{x}^{n+1}$ does not satisfy the following:

- The objective function must decrease or increase by, at most, 20%.
- The last Newton residual norm must decrease or increase by, at most, 20%.

$$F_V(\mathbf{x}^{n+1}) < 1.2 F_V(\mathbf{x}^n), \quad \left| \frac{\partial F_V(\mathbf{x}^{n+1})}{\partial \mathbf{x}_i} \right|_2 < 1.2 \left| \frac{\partial F_V(\mathbf{x}^n)}{\partial \mathbf{x}_i} \right|_2.$$

These modifications increase the probability of convergence significantly, i.e., the robustness of the method, but they may lead to more Newton iterations in the beginning stages of the nonlinear solve. The Newton method terminates when a positive $\alpha$ that satisfies the above requirements cannot be found.

**4. Additional practical aspects.** In this section we consider additional practical concerns beyond the quality of the optimized mesh.

**4.1. Tangential relaxation.** A common constraint on the mesh improvement process is that nodes belonging to the domain boundary should respect the boundaries of the domain (and potentially some prescribed internal curves/surfaces). The simplest way to accomplish this is to keep the boundary node coordinates *fixed* (i.e., unchanged) during the improvement process. At times, keeping boundary nodes fixed can be an annoying restriction because it unnecessarily limits the quality improvement that can be achieved. An alternative is to keep the boundary nodes on the domain boundary, but allow most of them to move tangentially to the boundary.

In tangential relaxation we still optimize the objective function from (2.3) or (2.4), but optimization of boundary nodes is treated differently. In 3D, surface nodes are allowed only two degrees of freedom, and curve nodes are allowed only one, whereas interior nodes have three. When the CAD description of the boundary is available, both derivative-free and derivative-based solvers optimize the objective functions of the boundary nodes in terms of their parametric coordinates. The positions of the mesh nodes are updated through the CAD parametrization. Derivative-based solvers require derivatives of the surface parametrization which are provided by the CAD model. Examples of such methods can be found in [35, 39].

The above approach is limited to the case in which the boundary curves are given parametrically. For cases in which the boundary is represented implicitly, one can still develop a similar tangential relaxation method in which the boundary nodes have less degrees of freedom, see [17].

Note that keeping the boundary nodes on the known boundary curve generally does not guarantee that the mesh boundary agrees with the curve, as the intermediate mesh positions depend on the chosen finite element space. In this paper we allow boundary motion only on boundaries that are parallel to the coordinate axes. In this case the motion restrictions are imposed as Dirichlet boundary conditions for the displacement vectors.

**4.2. Combination of metrics.** Using any of the existing metrics, we can define a composite metric as

$$\mu_{combo} = \eta_1 \mu_{i1}(T_{i1}) + \eta_2 \mu_{i2}(T_{i2}) + \dots \tag{4.1}$$

Each of the utilized metrics can have a different weight and different target matrix. The $\{\eta_k\}$ weights can be constant throughout the domain or a function of the location of the sample points. In time-dependent applications, they can also depend on time. By using a combination of metrics, different kinds of features such as shape, size and alignment can be emphasized in different regions of the mesh.

**4.3. Restricting the amount of mesh displacement.** In many applications with moving meshes (e.g. Lagrangian and ALE simulations) there is a trade-off between mesh quality improvement and accuracy of the field data transfer step (known as *remap*), see [4, 6]. While excessive node displacements may sometimes improve mesh quality, they might also impact the field transfer accuracy. The reason for this is that remap sub-steps, which are used to transfer a discrete field from one mesh configuration to another, introduce some amount of numerical diffusion (especially when the solution field is not smooth [3, 2]). Larger node displacements result in higher number of remap sub-steps, leading to increased diffusion.

Let $\mathbf{x}_k^{init}$ be the corresponding nodal coordinates of the initial mesh. Let $d \geq 0$ be some user-specified allowed displacement distance. Then, a natural way to formulate the restricted node movement problem would be to use the above objective functions, with the additional constraint

$$\left| \mathbf{x}_k - \mathbf{x}_k^{init} \right| \leq d \,.$$

Unfortunately, the numerical solution of problems with inequality constraints requires sophisticated solvers which can be expensive in practice. Here we present two alternative methods for restricting node movement during mesh optimization.

**4.3.1. First restricted node movement algorithm.** Let $\delta \geq 0$ be a user-supplied weighting parameter. Starting with the (2.4) objective function, we can define an alternative composite objective function of the form

$$\tilde{F_V}(x) = \sum_{E \in \mathcal{E}} \int_{E_t} \delta \left| x - x^{init} \right|^2 + \mu(T(x)) \, dx_t \,. \tag{4.2}$$

The optimization in the second restricted node movement algorithm is unconstrained (except for the usual fixed boundary nodes) and uses standard termination criterion. The movement away from the initial mesh will be reduced when $\delta$ is increased. A drawback of this approach is that the appropriate value of $\delta$ depends on the magnitude of $\mu(T)$ relative to the displacement distance.

**4.3.2. Second restricted node movement algorithm.** In this algorithm, node movement is restricted by using the following composite metric:

$$\mu'(T_0, T) = (1 - \gamma) \, \mu_{98}(T_0) + \gamma \, \mu(T) \,,$$

where

$$\mu_{98}(T) \equiv \frac{|T - I|^2}{\tau}$$

is the *replication* barrier metric.[2] Here, $T_0 = A W_0^{-1}$ and the target-matrix $W_0$ equals $A_{init}$, where $A_{init}$ is the Jacobian matrix at each sample point of the initial mesh.

---

[2] The metric *type* of $\mu_{98}$ is shape, size, and alignment.

Further, $T = AW^{-1}$, where $W$ is a different target-matrix which represents the desired Jacobian matrix in the unrestricted node movement problem (for example, $W = W_1$ could be the matrix corresponding to an ideal element shape). The metric $\mu'$ is thus really a function of $A$ and uses *two* target matrices, as well as two quality metrics. If $\gamma = 0$, then $\mu' = \mu_{98}$. At the first iteration of the optimization procedure, $\mu_{98} = 0$, i.e., the initial mesh is the minimizing mesh. Thus, optimization with $\gamma = 0$ replicates the initial mesh and no node movement occurs. On the other hand, if $\gamma = 1$, $\mu' = \mu$, and the optimization problem reduces to the original unrestricted node movement optimization problem. For intermediate values of $\gamma$, the optimization problem in the second algorithm is a mixture of the unrestricted and fully restricted problems.

**5. Numerical results.** In this section we report some results from the algorithms in the previous sections as implemented in the MFEM finite element library [31]. This implementation is freely available at `http://mfem.org`.

All results below are calculated using the *variational* objective function (2.4) and Gauss-Lobatto quadrature for the resulting integrals. Newton's method, as described in Section 3.4, is used to solve the nonlinear optimization problems. The Newton relative tolerance is set to $10^{-12}$, and all of the presented results are fully converged to this tolerance. The linear solve inside each Newton step is performed by the standard minimum residual (MINRES) algorithm. Boundary nodes which are on curved boundaries are always fixed, while the rest of the nodes are allowed to move as long as the motion does not perturb the initial domain.

One of the features of the MFEM library is that automatically provides efficient MPI parallelization of the evaluations of functionals, such as (2.4), which are assembled element-by-element. Furthermore, MFEM includes algorithms for the efficient assembly of the Hessian of $F_V$ and for its preconditioning in parallel, see [11, 10, 31].

**5.1. High-order mesh for a turbine blade.** We start with the test case that was mentioned in Section 3.3, namely, a $Q_4$ turbine blade mesh that is optimized with respect to shape. Optimization is performed with respect to the shape metric $\mu_2$, and the isotropic square target matrix $W_1$. The integral computations utilize 6 quadrature points in each direction. The initial mesh, along with its metric values $\mu_2(T(x))$ are shown on Figure 5.1. Note that all initial internal edges in the mesh are straight by construction. The goal of the optimization is to derive appropriate node displacements and curvature, so that the resulting node positions minimize the integral (2.4).

The optimized mesh and its metric values are shown in Figure 5.2. The shape is improved considerably, with larger final metric values around the 2 vertices that have 5 neighboring elements. Tangential relaxation allows node-movement on boundaries parallel to the coordinate axes. The objective function is reduced by approximately 61%, from $F_V(x_0) = 170.792$ to $F_V(x) = 66.4977$. The maximum metric values in the initial and optimized meshes are 21.5 and 16.7, respectively.

In addition to shape optimization, a common objective in mesh optimization is to preserve some boundary layer. As the original mesh has a well-pronounced boundary layer, we can preserve it by restricting the mesh displacement, as in Section 4.3. We repeat the above test by setting $\delta = 5000$ in (4.2) to obtain the result in Figure 5.3. The objective function is reduced by approximately 26%, from $F_V(x_0) = 170.792$ to $F_V(x) = 126.28$. The maximum metric values in the initial and optimized meshes are 21.5 and 6.12, respectively.

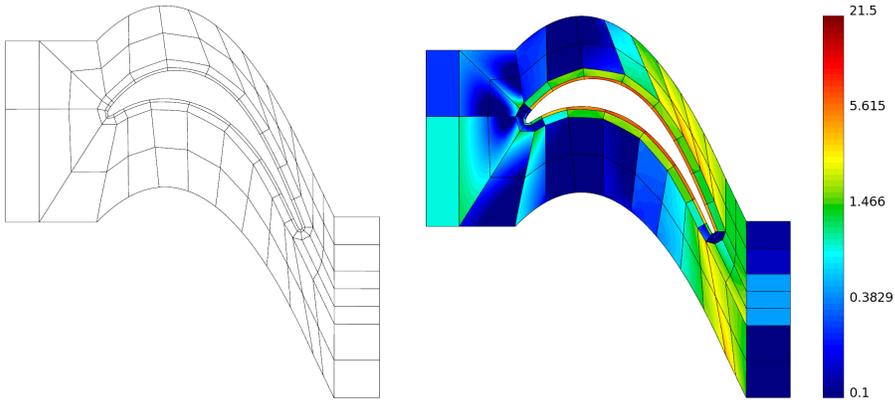

Fig. 5.1. *Initial $Q_4$ blade mesh and the corresponding $(\mu_2, W_1)$ metric values (log scale).*

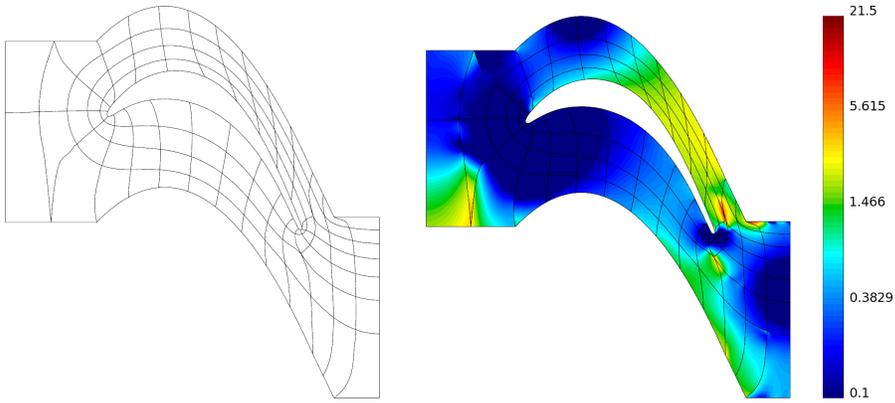

Fig. 5.2. *Optimized $Q_4$ blade mesh and the corresponding $(\mu_2, W_1)$ metric values (log scale).*

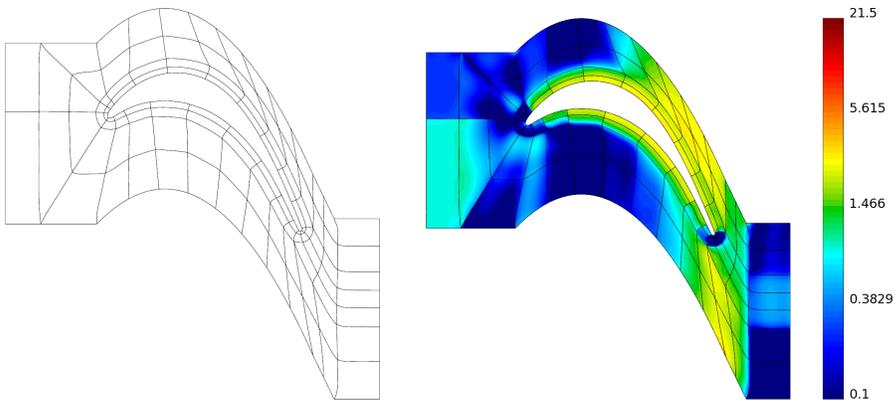

Fig. 5.3. *Optimized $Q_4$ blade mesh, with restricted mesh displacements, and the corresponding $(\mu_2, W_1)$ metric values (log scale).*

**5.2. High-order ICF meshes.** The following tests use a perturbed $Q_3$ mesh, see Figure 5.4. Similar meshes are often used for high-order ALE simulations [6] of Inertial Confinement Fusion (ICF) experiments, where preservation of radial symmetry is important. We use the following tests to demonstrate various capabilities of the method, e.g., different target constructions for shape and size optimization, composition of metrics, and restricting the amount of mesh displacement.

**5.2.1. Shape + equal size optimization.** In this test we use the $SS$ metric $\mu_9$ combined with the $W_2$ target, i.e., the mesh is optimized with respect to shape, so that the local size is uniform throughout the mesh. The integral computations utilize 6 quadrature points in each direction. The initial mesh, along with its metric values $\mu_9(T(x))$ are shown on Figure 5.4.

The optimized mesh and its metric values are shown in Figure 5.5. We observe good shape and symmetry, while the larger final metric values are at positions that are slightly above or below the average local size. The objective function is reduced by approximately 70%, from $F_V(x_0) = 1.119$ to $F_V(x) = 0.329$. The maximum metric values in the initial and optimized meshes are 298.6 and 28.58, respectively.

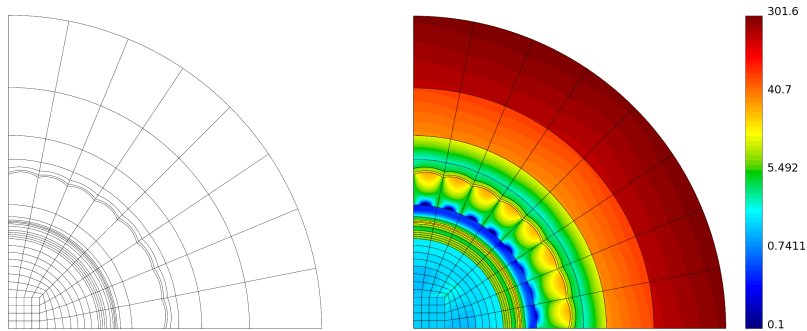

Fig. 5.4. *Initial $Q_3$ ICF mesh and the corresponding $(\mu_9, W_2)$ metric values (log scale).*

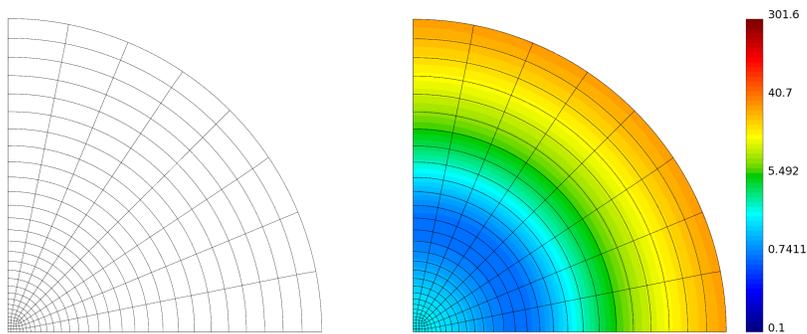

Fig. 5.5. *Optimized $Q_3$ ICF mesh and the corresponding $(\mu_9, W_2)$ metric values (log scale).*

**5.2.2. Shape + initial size optimization.** In this test we use the $SS$ metric $\mu_9$ combined with the $W_3$ target, i.e., the mesh is optimized with respect to shape, so that the local size around a given mesh node equals the local size around the node's initial position. The integral computations utilize 6 quadrature points in each direction.

The initial metric values, optimized mesh, and its metric values are shown in Figure 5.6. We observe improved shape and symmetry, while the initial local sizes are mostly preserved, except near the lower left corner. The objective function is reduced by approximately 51%, from $F_V(x_0) = 0.095$ to $F_V(x) = 0.046$. The maximum metric values in the initial and optimized meshes are 50.19 and 31.37, respectively.

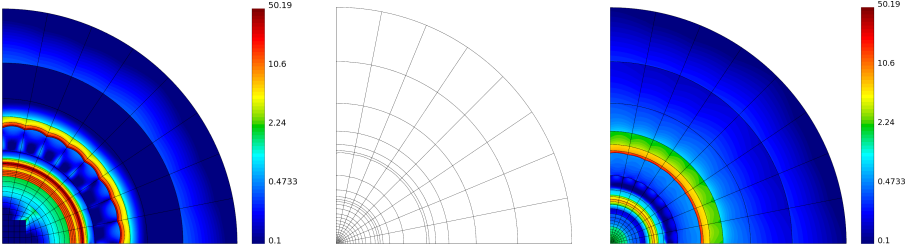

Fig. 5.6. *Initial metric values (left), optimized $Q_3$ ICF mesh (center), and the corresponding $(\mu_9, W_3)$ metric values (log scale, right).*

**5.2.3. Shape and restricted node movement optimization.** Next, we optimize the shape of the $Q_3$ ICF mesh while also trying to keep the mesh displacement small, as discussed in Section 4.3. We use the metric $\mu_1$ combined with the $W_1$ target. To reduce the mesh displacement we utilize the (4.2) objective function. The integral computations utilize 6 quadrature points in each direction.

To stress the effects of the restriction, results with $\delta = 0$ are presented, i.e., the mesh is optimized with no considerations about the node movement. The result is shown in Figure 5.7, which contains the initial $\mu_1$ metric values, the optimized mesh, and its $\mu_1$ metric values. The shape is improved, but all elements are substantially displaced and resized. The objective function is reduced by approximately 25%, from $F_V(x_0) = 0.189$ to $F_V(x) = 0.141$. The maximum metric values in the initial and optimized meshes are 0.0083 and 0.0070, respectively.

Results obtained with $\delta = 10$ are shown in Figure 5.8. The shape of the problematic elements is improved, while they remain close to their original positions. The position displacement (left panel) is localized around the regions with bad initial shape quality. The objective function is reduced by approximately 2.76%, from $\tilde{F}_V(x_0) = 0.189$ to $\tilde{F}_V(x) = 0.184$.

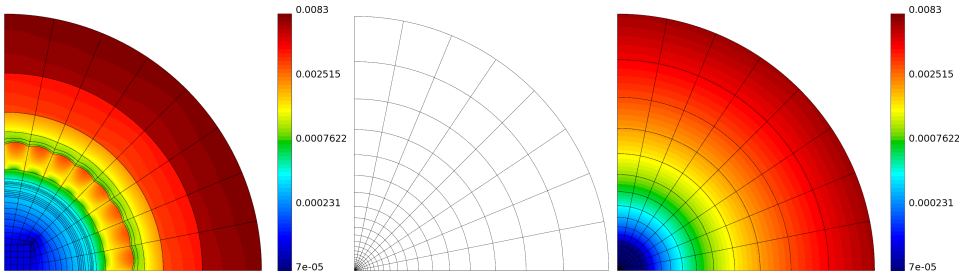

Fig. 5.7. *Initial metric values (left), optimized $Q_3$ ICF mesh (center), and the corresponding $(\mu_1, W_1)$ metric values (log scale), $\delta = 0$ (right).*

**5.2.4. Combination of metrics.** The goal of the next test is to optimize shape, while keeping higher resolution in a particular region of the domain. The following

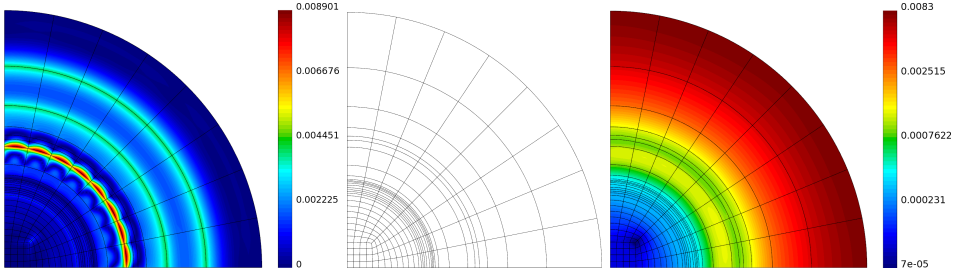

Fig. 5.8. *Amount of mesh displacement between the initial and optimized meshes (left), optimized $Q_3$ ICF mesh (center), and the corresponding $(\mu_1, W_1)$ metric values (log scale), $\delta = 10$ (right).*

objective function is used to optimize the perturbed ICF $Q_3$ mesh:

$$F_V(x) = \sum_{E \in \mathcal{E}} \left[ \int_{E_t} \mu_1(x) dx_t + \int_{E_t} \eta(x) \mu_{77}(x) dx_t \right],$$

where the space-dependent weight function $\eta$ is defined in terms of radius:

$$\eta(r) = 0.2 + 0.5 \left[ \tanh\left(\frac{r - 0.16}{0.002}\right) - \tanh\left(\frac{r - 0.17}{0.002}\right) + \tanh\left(\frac{r - 0.23}{0.002}\right) - \tanh\left(\frac{r - 0.24}{0.002}\right) \right].$$

The $\mu_1$ part of the objective function is used to optimize the mesh with respect to shape. The metric is combined with the $W_1$ target construction. The $\mu_{77}$ contribution optimizes the mesh with respect to size. It uses a target matrix of the form $W_2 = \sqrt{0.01\bar{\alpha}}\, I$, see (3.1). Such metric aims to achieve equal local sizes that are smaller than the average local size of the mesh. As the space-dependent weight $\eta(r)$ for $\mu_{77}$ stresses particular regions, the optimization tries to achieve finer resolution in these regions. This combination effectively forces shape and equal size optimization in most of the domain, while the stressed (by $\eta(r)$) regions are optimized with respect to shape and smaller size. The integral computations utilize 6 quadrature points in each direction.

The $\eta$ weight and final mesh are shown in Figure 5.9. The shape is improved throughout the mesh, while the regions, specified by $\eta(r)$, are resolved better. The final $\mu_1$ and $\mu_{77}$ metric values are shown in Figure 5.10. As $\mu_1$ has a constant weight in space, its final values are mostly uniform throughout the domain. The plot of $\mu_{77}$ shows that the stressed regions are optimized better in terms of the desired size. The objective function is reduced by approximately 80%, from $F_V(x_0) = 7.291$ to $F_V(x) = 1.421$.

**5.3. 3D mesh for a pinched sphere.** The following tests use 3D perturbed $Q_4$ and $Q_2$ spheres, see Figure 5.11. The internal shape of the meshes is optimized with metric $\mu_{303}$ and the $W_1$ target construction. This test demonstrates the benefits of global mesh optimization and the method's 3D mesh optimization functionality.

The $Q_4$ mesh (448 elements) results are shown in the top panel of Figure 5.11. The integral computations utilize 6 quadrature points in each direction. The shape of the inside of the sphere is substantially improved. The objective function is reduced by approximately 82%, from $F_V(x_0) = 169.605$ to $F_V(x) = 29.263$.

The $Q_2$ mesh (229,376 elements) results are shown in the bottom panel of Figure 5.11. The integral computations utilize 4 quadrature points in each direction. To verify the parallel execution capabilities of our MFEM implementation, a finer $Q_2$ mesh (1.835 million elements) was also tested on 512 MPI tasks, producing similar results (not shown here due to the density of elements).

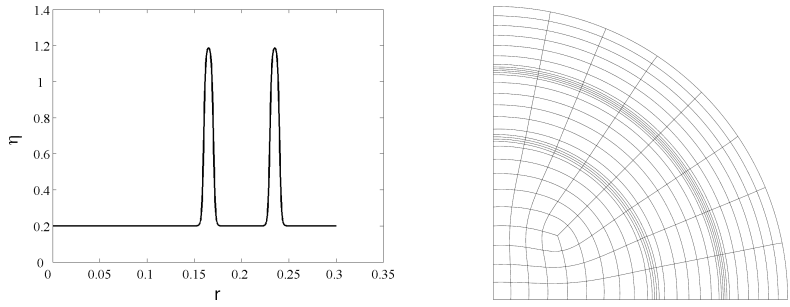

Fig. 5.9. *Space-dependent weight applied to* $\mu_{77}$ *(left) and optimized mesh (right).*

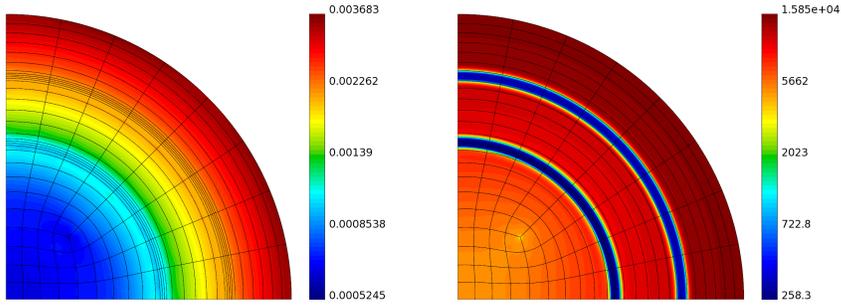

Fig. 5.10. *Final* $(\mu_1, W_1)$ *and* $(\mu_{77}, W_2)$ *metric values (log scale).*

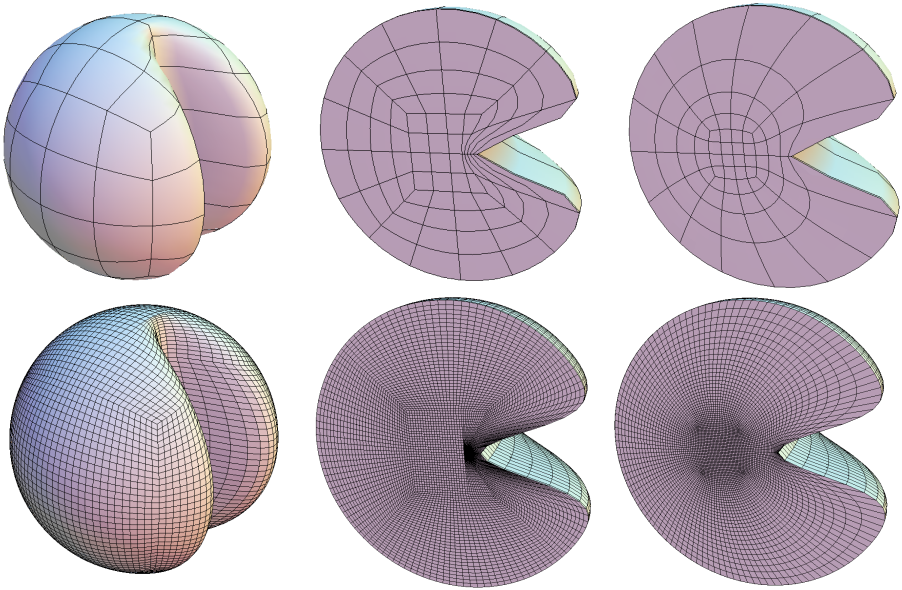

Fig. 5.11. *3D* $Q_4$ *(top) and a* $Q_2$ *(bottom) perturbed and optimized meshes (slices are shown in the middle and right panels).*

**6. Conclusions and future work.** This paper helps fill a gap in the high-order mesh literature on mesh quality improvement. This has been accomplished using the Target-matrix optimization paradigm via the important step of including sub-element information that is used to define the local quality of high-order elements. The paper demonstrates that this extension is sufficient to provide a degree of control over the quality of high-order 2D and 3D meshes in terms of local shape and area. The numerical method associated with high-order TMOP is implemented within the MFEM library [31] to provide the wider community free access to this method. Our approach is flexible and easy to incorporate in a wide range of current and future high-order applications, as well as low-order applications that use high-order meshes.

Although this method is already being used in realistic applications such as the BLAST code [6], considerably more work is needed to fully demonstrate the method, particularly in 3D and on other realistic applications. In principle, the target-matrices in TMOP provide the means to not only control local shape and size, but also to adapt the mesh to the physical solution as it evolves within the simulation (r-adaptivity). Some preliminary steps towards the goal of high-order r-adaptivity are reported in [15]; we hope to do much more regarding this challenge in the future. Other topics that will be addressed in our future work include: (1) tangential relaxation along curved interfaces and domain boundaries which are not given analytically, (2) application to realistic 3D simulations, (3) improvements in robustness and efficiency of the numerical algorithm and solvers, (4) investigation of additional target construction methods, (5) testing on applications which can benefit from the use of shape+size+alignment (SSA) metrics, and (6) coupling between node movement and adaptive mesh refinement (AMR).